\documentclass[12pt]{amsart}
\usepackage{amscd,amssymb}
\newtheorem{thm}{Theorem}[section]
\newtheorem{corol}[thm]{Corollary} 
\newtheorem{lemma}[thm]{Lemma}
\newtheorem{prop}[thm]{Proposition}
\theoremstyle{definition}
\newtheorem{defin}[thm]{Definition}

\theoremstyle{remark}

\numberwithin{equation}{section}

\def\norm#1{\left\Vert#1\right\Vert}
\def\norm#1{\left\Vert#1\right\Vert}
\def\oskip{\par\vbox to4mm{}\par}
\def\R{{\Bbb R}}
\def\Z{{\Bbb Z}}

\def\N{{\Bbb N}}
\def\T{{\Bbb T}}
\begin{document}
\title{SUBGROUPS OF MONOTHETIC GROUPS}

\author[S.A. Morris]{Sidney A. Morris}
\address{S.A.M.: School of Mathematics,
University of South Australia,
G.P.O. 2471, Adelaide, S.A., 5001, Australia}
\email{Sid.Morris@unisa.edu.au}
\date{\today}

\author[V. Pestov]{Vladimir Pestov}
\address{V.P.: School of Mathematical and Computing Sciences,
Victoria University of Wellington, P.O. Box 600, Wellington, 
New Zealand}
\curraddr{Computer Sciences Laboratory, RSISE, Australian National University,
Canberra, ACT 0200, Australia (July 1, 1999 -- Feb. 1, 2000)}
\email{vova@mcs.vuw.ac.nz}

\keywords{Monothetic groups, free abelian topological groups,
free locally convex spaces, Tkachenko--Uspenskij Theorem, 
$\omega$-tori, Rolewicz's Lemma}

\subjclass{22A05}

\begin{abstract}
It is shown that every separable abelian topological
group is isomorphic with a topological subgroup of a monothetic group
(that is, a topological group with a single topological generator).
In particular, every separable metrizable abelian group embeds
into a metrizable monothetic group.
More generally, we describe all topological
groups that can be embedded into monothetic groups: they are
exactly abelian
topological groups of weight $\leq\frak c$ covered by
countably many translations of every nonempty open subset.
\end{abstract}

\maketitle

\section{Introduction} 
A recent result by the present authors \cite{MP} states that
every separable topological group is isomorphic with a topological
subgroup of a group topologically generated by two elements.
This is a topological analogue of the Higman--Neumann--Neumann
Theorem. It leads to a simple description of those topological groups
which are embeddable into topological groups with two 
(equivalently, finitely
many) generators: they are exactly the topological groups that
are $\aleph_0$-bounded (that is, covered by countably many translations
of each non-empty open subset) and have weight $\leq\frak c$. 

Topological groups having one generator rather than two, that
is, monothetic groups (local compactness is not assumed here), 
are, naturally, abelian, and so are all their topological subgroups.
Rather surprisingly, abelianness turns out to be the only restriction
imposed on the previous result
in reducing the number of generators from two to one. 

\begin{thm}
\label{main}
Every separable abelian topological group
embeds into a singly generated topological group.
 \qed
\end{thm}

It is interesting to notice that, unlike the above mentioned topological
version of the Higman--Neumann--Neumann
theorem, our Theorem \ref{main} has no apparent discrete
algebraic counterpart.

\begin{corol}
\label{subgroups}
A topological group $G$ is isomorphic with a topological subgroup
of a monothetic group if and only if $G$ is abelian,
$\aleph_0$-bounded, and has weight $\leq\frak c$. \qed
\end{corol}

The present article was largely stimulated by the following, 
very general, question by Mycielski \cite{My}: {\it what can be said about
completely metrizable monothetic groups?} The following result 
gives some idea of the large size of such groups.

\begin{thm}
\label{metr}
Every separable metrizable abelian topological group is
isomorphic with a topological subgroup of a completely metrizable
monothetic group. 
\qed
\end{thm}

For example, all additive topological groups of separable Banach spaces
are to be found among subgroups of complete metric monothetic groups.

\oskip\section{Preliminary results and constructions}

Let $X=(X,\ast)$ be a pointed set, that is, a set with a
distinguished element $\ast\in X$. By $A(X,\ast)$, or simply
$A(X)$ if no confusion can arise, we shall 
denote the free abelian group on $X\setminus\{\ast\}$, 
having $\ast$ as its zero element,
and by $L(X,\ast)$ a real vector space having
$X\setminus\{\ast\}$ as its Hamel basis and $\ast$ as its zero
vector. It is well known and easily verified that
$A(X,\ast)$ is canonically isomorphic with a subgroup of the 
additive group of $L(X,\ast)$, generated by $X$. 

Let $\rho$ be a pseudometric on $X$. As was shown by Graev 
\cite{Gr}, there exists a maximal translation-invariant pseudometric
$\bar\rho$ on $A(X,\ast)$, whose restriction to $X$ coincides with
$\rho$.
A similar result was established \cite{AE, Rai} for
the vector span of $X$: there exists a maximal seminorm, $p_\rho$, on
$L(X,\ast)$, such that for every $x,y\in X$ one has
$\rho(x,y)=p_\rho(x-y)$.

Since the pseudometric on $A(X)$ induced by the seminorm $p_\rho$ is
clearly translation invariant, for all
$x,y\in A(X,\ast)$ one has $p_\rho(x-y)\leq\rho(x,y)$.
The following important and nontrivial
result, obtained by successive efforts of Tkachenko
\cite{T} and
Uspenskij\footnote{Other transliterations used: Uspenski\u\i\ and Uspenskiy.} 
(\cite{U}, pp. 660--662), shows
that the two pseudometrics on $A(X)$ thus obtained in fact coincide.

\begin{thm}[Tkachenko--Uspenskij] If  $\rho$ is any pseudometric
on a pointed set $(X,\ast)$, then
$p_\rho(x-y)=\bar\rho(x,y)$ for all $x,y\in A(X,\ast)$.
\qed
\label{usp}
\end{thm}

If now $X$ is a Tychonoff topological space, then the
{\it free abelian topological group} on $X$ is the group $A(X)$
equipped with the finest group topology inducing the given topology
on $X$ as a subspace. Such a topology always exists, is Hausdorff,
and has the universal property of the following kind: every
continuous mapping $f$ from $X$ to an arbitrary abelian topological
group $G$ lifts to a unique continuous homomorphism
$\bar f\colon A(X)\to G$. It was first observed by Graev that
the topology of $A(X)$ is determined by the collection of all
translation invariant pseudometrics of the form $\bar \rho$, where
$\rho$ is a continuous pseudometric on $X$. For an account of theory
of free abelian topological groups, see e.g. \cite{Mo}.

In a similar way, the {\it free locally convex space} on $X$ is
the vector space $L(X)$ equipped with the finest locally convex
topology inducing the given topology on $X$. Such a topology exists
and is Hausdorff whenever $X$ is a Tychonoff topological space,
and every continuous mapping $f$ from $X$ to an
arbitrary locally convex space $E$ extends to a unique continuous
linear operator from $L(X)$ to $E$. The topology of $L(X)$ is
determined by the collection of all seminorms of the type $p_\rho$ 
\cite{AE, Rai}.

The following result was obtained by Graev \cite{Gr}.

\begin{thm}[Graev] Let $G$ be a Hausdorff topological group,
and let $i\colon A(G)\to G$ be the unique continuous homomorphism
from the free abelian topological group on (the underlying
topological space of) $G$ to the topological group $G$ whose
restriction to $G$ is the identity map. Then $i$ is a quotient homomorphism
of topological groups.
\qed
\label{graev}
\end{thm}

The Tkachenko--Uspenskij Theorem was put to use in \cite{MM}
where the following technique was suggested. Let $G$ be a
topological group. Denote by $i\colon A(G)\to G$ the homomorphism
described in Graev's Theorem \ref{graev}, and let $K$ denote the
kernel of $i$. Then $K$ is a closed topological subgroup of $A(G)$,
and the topological factor-group $A(G)/K$ is isomorphic to $G$.
Moreover, the same remains true if we consider the group
$A(G)$ equipped with the Graev metric $\bar d$, where $d$ is
a translation-invariant metric generating the topology of $G$.
According to Tkachenko--Uspenskij Theorem,
$(A(G),\bar d)$ is isomorphic to a topological subgroup of
the normed space $(L(G),p_d)$.
It is now easy to see that the topological factor-group $L(G)/K$
contains $G$ as a (closed) topological subgroup. 

\begin{prop}[\cite{MM}] Every metrizable abelian topological group
is isomorphic to a topological subgroup of a topological
factor-group of the additive group of a suitable
Banach space.
\qed
\end{prop}

We need to develop a slight technical variation on the above themes.

Let $\frak R$ be a collection of pseudometrics on a pointed set
$X=(X,\ast)$ inducing some (Tychonoff) topology. (In precise terms, the
collection of all open balls formed with respect to pseudometrics
from $\frak R$ forms a topology base.) We will denote by
$A(X,\ast, {\frak R})$ the free abelian group $A(X,\ast)$, equipped with
the collection of Graev pseudometrics 
$\bar{\frak R}=\{\bar\rho\colon \rho\in {\frak R}\}$.
Mostly we shall be viewing $A(X, \ast,{\frak R})$
as an abelian topological  group under the (group) topology generated
by all pseudometrics from $\bar{\frak R}$. 
The space $X$ evidently is a topological subspace of
$A(X, \ast,{\frak R})$.
The following is also clear. Recall that a mapping $f\colon X\to Y$
between two metric spaces $X=(X,\rho_X)$ and $Y=(Y,\rho_Y)$ is
called {\it Lipschitz} if for some $C\geq 0$ one has
$\rho_Y(f(x),f(y))\leq C \rho_X(x,y)$ for all $x,y\in X$.
Such a $C$ is called a {\it Lipschitz constant} for the mapping $f$.
If $f$ is Lipschitz with Lipschitz constant $C=1$, then $f$ is
called a {\it 1-Lipschitz mapping.}

\begin{prop} Let $(X,\ast, {\frak R})$ be as above.
Let $G$ be an abelian topological  group, and let
$\frak P$ be some collection of translation-invariant pseudometrics
on $G$
generating the topology.
Let $f\colon X\to G$ be a mapping sending $\ast$ to $0_G$ and such
that for every $d\in {\frak P}$ there is a $\rho\in{\frak R}$
making the map $f\colon (X,\rho)\to (G,d)$ Lipschitz.
Then the (unique) algebraic homomorphism $\bar f\colon A(X)\to G$,
extending $f$, is continuous. \qed
\label{unia}
\end{prop}

In a similar way, we define a locally convex space
$L(X,\ast, {\frak R})$ as the linear span $L(X,\ast)$
of $X$ with $\ast$ serving as zero equipped
with the (locally convex topology generated by) the collection of
seminorms $\widetilde{\frak R}=\{p_\rho\colon \rho\in {\frak R}\}$.
The space $X$ is a topological subspace of
 $L(X,\ast, {\frak R})$. 
One has the following counterpart of Proposition \ref{unia}.

\begin{prop} Let $E$ be a locally convex space, and let
$\frak P$ be some collection of seminorms
generating the topology of $E$.
Let $f\colon X\to E$ be a mapping sending $\ast$ to $0_E$ and such
that for every $p\in {\frak P}$ there is a $\rho\in{\frak R}$
making the map $f\colon (X,\rho)\to (E,p)$ Lipschitz.
Then the (unique) linear operator 
$\bar f\colon L(X,\ast, {\frak R})\to E$,
extending $f$, is continuous. \qed
\label{unil}
\end{prop}

Now
Tkachenko--Uspenskij Theorem \ref{usp} implies the following.

\begin{corol}
\label{new} The additive topological group of
the locally convex space $L(X,\ast,{\frak R})$ contains an
isomorphic copy of the abelian topological group
$A(X, \ast, {\frak R})$ in a canonical way.
\qed
\end{corol}

Moreover, the copy of $A(X)$ forms a closed topological subgroup of $L(X)$
({\it ibid.})

If $X=(X,\ast)$ is a pointed Tychonoff space
(that is, a Hausdorff topological space in which points
and closed subsets are separated by continuous real-valued
functions), we will
denote by ${\frak R}(X)$ the collection of all continuous 
pseudometrics on $X$.
Then $A(X,\ast,{\frak R}(X))$ is naturally isomorphic to the Graev free
abelian topological group, $A(X)$, on $X$, cf. \cite{Gr, Mo}, while
$L(X,\ast,{\frak R}(X))$ is naturally isomorphic to the Graev free
locally convex space, $L(X)$, on $X$, cf. \cite{AE, Fl, Rai}.

Let $G$ be an abelian topological group, and let $\frak R$ be an
arbitrary
family of continuous translation-invariant 
preudometrics generating the topology of $G$.
The identity mapping $\operatorname{Id}_G$ satisfies the assumption
of Proposition \ref{unia} and therefore extends to a unique
continuous surjective homomorphism 
 $i\colon A(X,\ast,{\frak R})\to G$,
sending $\ast$ to $0$.

\begin{prop}
The homomorphism $i\colon A(X,\ast,{\frak R})\to G$ is open.
\end{prop}

\begin{proof}
In fact, the same homomorphism $i$ is open even if
considered as a mapping
from the free abelian topological group $A(G)$ to $G$ \cite{Arh}.
\end{proof}

Denote by $K_G$ the kernel of $i$, which is a closed topological
subgroup of $A(X,\ast,{\frak R})$. The openness of $i$ implies the
following.

\begin{corol}
$G$ is canonically topologically isomorphic to the topological factor
group $A(X,\ast,{\frak R})/K_G$.
\qed
\end{corol}

Let $G$ be a countable topological group. Now choose as
$\frak R$ the collection of all translation-invariant continuous
pseudometrics on $G$ which are bounded by $1$. Then $\frak R$ determines
the topology of $G$ (which is true of every topological group $G$,
cf. \cite{Gr}). Denote by $d$ the discrete metric on $G$, that is, one
taking values $0$ and $1$ only. (In general, $d$ is discontinuous ---
unless of course $G$ is discrete.) Notice that for each $\rho\in
{\frak R}$ the identity mapping $(G,d)\to (G,\rho)$ is 1-Lipschitz.
It implies that the identity isomorphism
$(A(G),\bar d)\to (A(G), {\frak R})$ is continuous.
Noticing that the locally convex space $L(G,d)$ is separable metrizable
and contains $(A(X),\bar d)$ as a topological subgroup
(Corollary \ref{new}),
one arrives at the following.

\begin{lemma}
Every separable topological group $G$ is isomorphic to a
topological factor-group of a group $A(G,{\frak R})$,
contained as a topological subgroup in a separable locally convex
space $L(G,{\frak R})$, admitting a finer separable metrizable
locally convex topology. \qed
\label{glav}
\end{lemma}

\oskip\section{Rolewicz's lemma}

Let us introduce the following convenient notion.

\begin{defin} We call an abelian topological group $G$ an
{\it $\omega$-torus} if it is topologically generated by
the union of countably infinitely many subgroups topologically
isomorphic to the circle group ${\Bbb T}\cong
U(1)$. 
\end{defin}

The following was essentially proved by Rolewicz
\cite{Ro}. Even though he did not  state the result in its full
generality, the proof is his. The construction forms a
rich source of monothetic groups beyond the locally compact case
(cf. e.g. \cite{DPS} and references therein). We therefore find it
very useful, to state Rolewicz's Lemma in its full generality,
and believe that such a generalization is of interest on its own
and not just in connection with the subsequent application in this
article.
In the proof we will stick to the multiplicative notation as more
convenient in this particular context.

\begin{thm}
\label{mon}
Every completely metrizable $\omega$-torus is monothetic. 
\end{thm}

\begin{proof}
Let $G$ be  topologically
generated by the union of a countable sequence of its subgroups
$\T_i$, $i=1,2,\dots$, each of which is topologically isomorphic to
the circle group $\T=U(1)$.
Fix a translation invariant metric, $\rho$, generating the
topology on $G$.
For each $i=1,2,\dots$ choose recursively
a number $n_i\in\N$
and an element $x_i\in {\Bbb T}_i$ satisfying the following properties
for each $i$.

\begin{enumerate}
\item $\rho(x_i,0)<2^{-i}$.
\item The first $n_i$ powers of the
 product 
$x_1\cdot x_2\cdot\dots\cdot x_i$
form a $2^{-i}$-net in the compact subgroup 
$\T_1\cdot \T_2\cdot\dots \cdot\T_i$ of $G$.
\item Whenever $j>i$, the first $n_i$ powers of the element $x_j$
are contained in the $\rho$-ball of zero
having radius $2^{-j}$.
\end{enumerate}

As the base of recursion, choose any element $x_1\in\T_1$ having infinite
order and contained in the $1/2$-neighbourhood of zero formed with respect
to the metric $\rho$. To perform the recursive step, 
recall the classical Kronecker Lemma: if $x_1,\dots,x_n$ are
rationally independent real numbers, then the $n$-tuple 
$(x_1',\dots,x_n')$ made up of their images under the quotient homomorphism
$\R\to\R/\Z$ to the circle group $\T\cong\R/\Z$ generates an
everywhere dense subgroup in the $n$-torus ${\Bbb T}^n$.
Now assume that $x_1,\dots,x_{i-1}$ and $n_1,\dots,n_{i-1}$ 
with the properties
(1)--(3) have been chosen. If the closed subgroup, $A_i$, 
generated by the product
$x_1\cdot\dots\cdot x_{i-1}$ coincides with all of $T_1\cdot\dots\cdot T_i$,
we set $x_i:=0$. Otherwise, $A_i$ forms a proper closed subgroup of the
group $T_1\cdot\dots\cdot T_i$, and clearly
the latter is isomorphic to the topological direct sum $A_i\times \T_i$.
Moreover, the compact abelian Lie group $A_i$ is itself isomorphic
to a torus group $\T^j$ of a suitable
rank $j\leq i-1$, and the image of the 
topological generator $x_1\cdot\dots\cdot x_{i-1}$ in $\T^j$
under such an isomorphism
is a $j$-tuple, say $(z_1,\dots,z_j)$, 
of rationally independent elements of the circle group.
Now choose $x_i\in \T_i$ to be an element that is rationally independent
of the elements $z_1,\dots,z_j$ and such that all the powers
$x_i,x_i^2,\dots,x_i^{n_{i-1}}$ are contained in the $\rho$-neighbourhood
of zero having radius $2^{-i}$. It follows from the Kronecker Lemma
that the powers of the product $x_1\cdot\dots\cdot x_i$ are everywhere
dense in the group $\T_1\cdot\dots\cdot \T_i$.
Consequently, it is possible to choose $n_i$ as a sufficiently
large natural number so that the first $n_i$ powers of
 $x_1\cdot\dots\cdot x_i$ form a $2^{-i}$-net in $\T_1\cdot\dots\cdot \T_i$.
The step of recursion is thus accomplished.

We claim that the element $x=\prod_{l=1}^\omega x_l$,
which is clearly well-defined since the metric $\rho$ on 
$G$ is complete, forms a topological generator for
the group $G$. It is enough to demonstrate that for every number
$k\in\N$ the closure of the cyclic group $\langle x\rangle$
generated by $x$ contains $\T_k$.
Let $k\in\N$ and let $z\in \T_k$ by any. 
Let now $i\geq k$ be arbitrary. Since
$z\in \T_1\cdot\T_2\cdots \T_i$, condition (2) implies the existence
of an $m=0,1,2,\dots,n_i$
such that the $m$-th power of
$x'=x_1\cdot x_2\cdot\dots\cdot x_i$ is at a distance $<2^{-i}$ from
$z$. The $m$-th power of the remainder of the infinite product, 
$x(x')^{-1}=\prod_{l=i+1}^\omega x_l$, is at a distance
from zero which is less than
\[\sum_{l=i+1}^\omega \rho(0,(x_l)^m)<\sum_{l=i+1}^\omega 2^{-l}=2^{-i}.\]
(Here we have used condition (3).) Finally,
\[\rho(x^m,z)\leq\rho((x')^m,z)+\rho((x(x')^{-1})^m,0)<2^{i-1}.\]
Since $i$ can be chosen arbitrarily large, 
the latter inequality means that $z$ is the limit of a sequence of
suitable powers of $x$, and the proof is finished.
\end{proof}

\oskip\section{The main construction}

Assume we are given the following collection of data.

\begin{enumerate}
\item A separable topological vector space $E$.
\item A countable everywhere dense subset of $E$, 
$X=\{x_m\colon m\in\N_+\}$.
\end{enumerate}

Form the direct sum topological vector space, equipped with the direct
product topology:
\[H:=E\oplus l_2(\N_+).\]
We will identify $E$ 
in a natural way with the topological vector subspace (and subgroup) 
first direct summand of $H$.
For each pair $(m,n)$ of positive 
natural numbers, denote 
\[\xi_{m,n}:=(nx_m,e_{m,n})\in H.\]
 Let $D$ denote the subgroup of $H$ algebraically generated by all
elements  $\xi_{m,n}, m,n\in\N_+$. 

Let $d$ be an arbitrary translation invariant continuous pseudometric
on $E$. Such pseudometrics generate the topology of $E$. Denote by
$\tilde d$ the continuous translation invariant pseudometric on $H$
defined  by setting for each 
$x_1,x_2\in E$ and $y_1,y_2\in l_2(\N_+^2)$
\[\tilde d((x_1,y_1),(x_2,y_2)):= d(x_1,x_2)+\norm{y_1-y_2}.\]
The collection of all pseudometrics of the form $\tilde d$ generates
the topology of $H$. 

\begin{lemma}
The neighbourhoods of zero, $W$, with the property $(D+W)\cap E= W\cap E$
form a neighbourhood basis in $H$.
\label{ba}
\end{lemma}

\begin{proof} Let $d$ be an arbitrary 
pseudometric on $E$ as above, and let
\[W:=\{z\in H\colon \tilde d(z,0)<1\}.\]
It is clearly enough to show that $W$ has the required property, that is,
if $x\in D\setminus \{0\}$ and $y\in E$, then $\tilde d(x,y)\geq 1$.

An arbitrary element, $x$, of $D$ is of the form 
\begin{equation}
\sum_{i=1}^s k_i\xi_{m_i,n_i}\equiv \left(\sum_{i=1}^s k_i n_i x_{m_i},
\sum_{i=1}^s k_ie_{m_i,n_i}\right)\in H.
\end{equation}
Let $x\neq 0$. Then
one can assume without loss in generality
that in the above expansion
all the integer coefficients $k_i\neq 0$, and that for different
$i,j$ the pairs $(m_i,n_i)\neq (m_j,n_j)$. 
Now let $y\in E$ be arbitrary. According to our earlier
convention, we will identify $y$ with
the element $(y,0)\in H$.
One has
\begin{eqnarray} 
\tilde d(x,y)&=& d\left(\sum_{i=1}^s k_i n_i x_{m_i},y\right)+
\norm{\sum_{i=1}^s k_ie_{m_i,n_i}} \nonumber \\
&\geq & \norm{\sum_{i=1}^s k_ie_{m_i,n_i}} \nonumber \\
&\geq&  1,
\end{eqnarray}
and the claim follows.
\end{proof}

\begin{lemma}
The group $D$ is discrete.
\end{lemma}

\begin{proof}
Indeed, the image of $D$ under the second
coordinate projection $H\to l_2(\N_+)$ (which is of course a homomorphism of
topological vector spaces) is a discrete subgroup of $l_2(\N_+)$,
formed by all linear combinations of the standard basic elements
with integer coefficients.
\end{proof}

\begin{lemma}
The linear span of $D$ is everywhere dense in $H$.
\label{ed}
\end{lemma}

\begin{proof}
Let $m\in\N_+$ be arbitrary. For every $n\in\N_+$, the linear span
of $D$ contains the element $(1/n)\xi_{m,n}=(x_m,(1/n)e_{m,n})$, and
for each continuous pseudometric $d$ on $E$
\begin{eqnarray}
\tilde d\left(\frac 1n\xi_{m,n},x_m\right) &=&
d(x_m,x_m)+\norm{0,\frac 1n e_{m,n}}\nonumber \\
&= & \frac 1 n.
\end{eqnarray}
Consequently, $x_m$ is in the closed linear span of $D$. Since the set
$\{x_m\colon m\in\N_+\}$ is everywhere dense in $E$, it follows that the
closed linear span of $D$ contains $E$. Further,
each element of the form
\[e_{m,n}=\xi_{m,n}-\frac 1n x_m,\]
where $m,n\in\N_+$,
is in the closed linear span of $D$ as well. But $\{e_{m,n}\colon m,n\in\N_+\}$
is an orthonormal basis for $l_2(\N_+)$. The claim is established.
\end{proof}

\begin{lemma}
The factor-group $H/D$ is an $\omega$-torus.
\label{infty}
\end{lemma}

\begin{proof}
According to Lemma \ref{ed}, the topological group 
$H$ is topologically generated by the union of countably many
one-parameter subgroups passing through the elements of the form $\xi_{m,n}$. 
Therefore, $H/D$ is topologically generated
by the union of images of all such one-parameter subgroups. 
But those images are tori,
and there are countably many of them.
\end{proof}

\begin{lemma}
The restriction of the quotient homomorphism $\pi\colon H\to H/D$
to $E$ is a topological group 
isomorphism between $E$ and its image
in $H/D$.
\label{emb}
\end{lemma}

\begin{proof}
Since $D\cap E=\{0\}$, the homomorphism $\pi\vert_E$ is in
fact an algebraic isomorphism, and clearly it is continuous.
It remains to prove that $\pi\vert_E$ is open on its image.
Let $V$ be an arbitrary neighbourhood of zero in $H$.
According to Lemma \ref{ba}, there is a neighbourhood of zero
$W$ contained in $V$ with the property $(D+W)\cap E= W\cap E$,
that is, $\pi(W\cap E)=\pi(W)\cap\pi(E)$. Consequently,
the interior of $\pi(W\cap E)$ in $\pi(E)$ is non-empty
(it contains the open set $\pi(W)\cap\pi(E)$ in the subspace
topology induced from $H/D$), and the proof is finished.
\end{proof}

\begin{lemma}
\label{also}
If $E$ is also metrizable, then $H/D$ is monothetic metrizable.
\end{lemma}

\begin{proof}
The metrizability of $H$ --- and therefore of $H/D$ --- is quite
obvious, and monotheticity of $H/D$ follows from Rolewicz's Lemma
and Lemma \ref{infty}.
\end{proof}

Combining together Lemmas \ref{also} and \ref{emb}, 
we obtain the following.

\begin{lemma}
\label{tvs}
 Let $E$ be a separable metrizable topological vector space.
Then $E$ embeds, as a topological group, into a monothetic metrizable
group $H$.
\qed
\end{lemma}

Our next task will be to obtain a similar result for a sufficiently wide
class of separable non-metrizable topological vector spaces.

\begin{lemma}
\label{comp}
Let the TVS $E$ admit a finer
Hausdorff topology $\frak T$ that makes
it into a separable metrizable topological vector space.
Then the group $H/D$ is monothetic.
\end{lemma}

\begin{proof} Denote by $F$ the underlying vector space of $E$ equipped
with the topology $\frak T$. The identity mapping $F\to E$
is a continuous linear operator, and as such, it extends over the
completions of the two topological
vector spaces in a unique way, giving rise
to a continuous homomorphism (in general, no longer an algebraic
isomorphism)
\[i\colon \widehat F\to\widehat E.\]
Choose a countable everywhere dense subset $X\subset F$, 
$X=\{x_m\colon m\in\N_+\}$. Then $X$ remains everywhere dense in
$E$ as well.
Now apply our construction to both spaces $\widehat E$ and 
$\widehat F$. To distinguish
between the emerging pairs of 
objects, we will be using subscripts $E$ and $F$,
respectively. Thus, $H_F=\widehat F\oplus l_2(\N_+)$, etc.
Obviously, the subgroups $D_E$ and $D_F$ coincide as abstract groups,
or, more precisely, $i\vert_{D_F}\colon D_F\to D_E$ is a topological
group isomorphism. Consequently, the homomorphism $i$ factors through
$D_F$ to give rise to a continuous homomorphism 
$j\colon H_F/D_F\to H_E/D_E$. Quite evidently, the image of $j$
forms an everywhere dense subgroup of $H_E/D_E$.
Since the former of the two groups is completely
metrizable, Rolewicz's Lemma coupled with
Lemma \ref{infty} imply that $H_F/D_F$ is a monothetic group.
But the image of a monothetic group under a continuous homomorphism
having dense image is again monothetic.
\end{proof}

The following is a direct consequence of Lemmas \ref{comp}
and \ref{emb}.

\begin{lemma}
Every topological vector space $E$ that admits a finer
Hausdorff topology $\frak T$ making it
into a separable metrizable topological vector space
embeds as a topological subgroup into a monothetic group.
\label{teh}
\qed
\end{lemma}

\section{Proofs of the main results}

Denote provisionally by $\mathcal G$ the class of all 
topological groups that embed, as topological subgroups, into
monothetic groups.

The following is straightforward.

\begin{lemma}
The class $\mathcal G$ is closed under passing to topological subgroups.
\label{subgroups1}
\end{lemma}

Our next result is less evident.

\begin{lemma}
The class $\mathcal G$ is closed under passing to topological
factor-groups.
\label{quotients}
\end{lemma}

\begin{proof}
Indeed, let $G\in {\mathcal G}$, that is, for some monothetic group
$H$, one has $G<H$. Let $F$ be a closed subgroup of $G$. Denote by
$F'$ the closure of $F$ in $H$; one has $F'\cap G=F$.
Then it is a standard result, repeatedly used in theory of varieties
of topological groups (cf. e.g. \cite{survey}), that
the topological factor-group $G/F$ is isomorphic to a
topological subgroup of the factor-group $H/F$ in a canonical way.
At the same time, the group $H/F$ is clearly monothetic.
\end{proof}

\subsection*{Proof of Theorem \ref{main}}
Follows from Lemmas \ref{glav}, \ref{subgroups1}, \ref{quotients},
and \ref{teh}.
\qed

\subsection*{Proof of Corollary \ref{subgroups}}
It suffices to apply a description of topological subgroups of
separable topological groups obtained in \cite{MP}: those are
exactly $\aleph_0$-bounded topological groups of weight $\leq\frak c$.
Both the statement and the proof
remain true if we add the word `abelian' throughout. \qed

\subsection*{Proof of Theorem \ref{metr}}
Every separable metrizable abelian topological group,
$G$, is
isomorphic to a topological factor-group of the free abelian group
equipped with the Graev metric, $A(G,\rho)$, where
$\rho$ is an arbitrary metric on $G$ generating the topology.
By Tkachenko--Uspenskij Theorem, $A(X,\rho)$ is isomorphic with
a closed topological sugroup of the free Banach space,
$B(X,\rho)$. Notice that $B(X,\rho)$ is separable as well.
According to Lemma \ref{tvs}, $B(X,\rho)$ is isomorphic with a
topological subgroup of a suitable monothetic metrizable group,
say $H$. An application of Lemma \ref{subgroups1} and Lemma
\ref{quotients} finishes the proof.
\qed

\section*{Acknowledgements} 
The second named author (V.P.) is grateful to the Mathematical Analysis
Research Group of the Centre for Industrial
and Applicable Mathematics (CIAM)
of the University of South Australia for support
and hospitality extended during his visit in July 1998.
The research of the same author was in part supported
by the Marsden Fund grant VUW703 of the Royal Society of New Zealand.

\end{document}